\documentclass[10pt,twoside]{amsart}
\usepackage{lmodern}
\usepackage[T1]{fontenc}
\usepackage{amsmath,amsfonts,amsthm,amssymb,amscd}
\usepackage{mathrsfs}
\usepackage{euscript}
\usepackage{latexsym}
\usepackage[pdftex]{graphicx}
\usepackage{enumitem}
\usepackage{bbm}
\usepackage{accents}

\usepackage{geometry}
\newtheorem{lemma}{Lemma}
\newtheorem{theorem}{Theorem}

\usepackage{hyperref}  
\hypersetup{colorlinks=true, linktoc=all,
    linkcolor=blue,
		citecolor=blue,
		urlcolor=cyan,
		bookmarksopen=true
		}
\title{Estimates of short exponential  sums with primes in major arcs}
\author{Zarullo Rakhmonov}
\address{A.Dzhuraev Institute of Mathematics,  National Academy of Sciences of Tajikistan}
\email{zarullo.rakhmonov@gmail.com}
\date{12.11.2026}
\begin{document}

\begin{abstract}
When solving a number of additive problems with almost equal summands, along with estimating short trigonometric sums with primes of the form
$$
S_k(\alpha;x,y)=\sum_{x-y<n\le x}\Lambda(n)e(\alpha n^k),
$$
on the minor arcs, estimates of these sums on the major arcs are also needed, except in a small neighborhood of their centers, as well as an asymptotic formula in a small neighborhood of the center of the major arcs.

In this paper, using the second moment of Dirichlet $L$-functions on the critical line, for $S_k(\alpha;x,y)$ on the major arcs $\mathfrak{M}(\mathscr{L}^b)$, $\tau=y^5x^{-2}\mathscr{L}^{-b_1}$, $\mathscr{L} =\ln xq$, except in a small neighborhood of their centers $|\alpha-\frac{a}{q}|>\left(2\pi k^2x^{k-2}y^2\right)^{-1}$, for $y\ge x^{1-\frac{1}{2k-1+\eta_k}}\mathscr{L}^{c_k}$,
$$
\eta_k=\frac{2}{4k-5+2\sqrt{(2k-2)(2k-3)}}, \qquad c_k=
 \frac{2A+22+\left(\frac{2\sqrt{2k-3}}{\sqrt{2k-2}}-1\right)b_1}{2\sqrt{(2k-2)(2k-3)}-(2k-3)},
$$
a nontrivial estimate of the form
$$
S_k(\alpha;x,y)\ll y\mathscr{L}^{-A},
$$
is obtained, where $A$, $b_1$, $b$ are arbitrary fixed positive numbers, and an asymptotic formula is proved in a small neighborhood of the centers of the major arcs.
\end{abstract}

\maketitle

\section{Introduction}

According to Dirichlet's theorem on the approximation of real numbers by rational numbers, every $\alpha$ in the interval $\left[-\frac1\tau, 1-\frac1\tau\right]$, can be represented in the form
\begin{equation*}
 \alpha =\frac{a}{q}+\lambda ,\quad (a,q)=1,\qquad 1\le q\le \tau , \qquad |\lambda |\le \frac{1}{q\tau }.
\end{equation*}
Let $\mathfrak{M}(P)$ denote those numbers $\alpha$ for which $q\le P$, and let $\mathfrak{m}(P)$ denote the remaining $\alpha$. $\mathfrak{M}(P)$ and $\mathfrak{m}(P)$ are called the major and minor arcs, respectively.

Vinogradov~I.M.~\cite{VinigradovIM-Izb-Trud,Vinogradov+Karatsuba-TrMIAN-1984-77} was the first to study short trigonometric sums with primes. For sums of the form
\begin{align*}
&S_k(\alpha ;x,y)  = \sum_{x-y<n\le x}\Lambda(n) e(\alpha n^k),
\end{align*}
for $k=1$, using his method of estimating sums with primes, he proved a nontrivial estimate on the minor arcs $\mathfrak{m}(\exp(c(\ln \ln x)^2))$, $\tau =x^{\frac{1}{3}}$ for $y>x^{\frac{2}{3}+\varepsilon}$, whose basis, along with the "Vinogradov sieve", for $k=1$ consists of estimates of short double trigonometric sums of the form
$$
J_k(\alpha;x,y,M,N)=\sum_{M<m\le 2M}a_m\sum_{U<n\le 2N \atop x-y<mn\le x}b_ne(\alpha (mn)^k),
$$
where $a_m$ and $b_n$ are arbitrary real functions, $|a_m|\le \tau^c(m)$, $|b_n|~\le~\tau^c(n)$, $M$, $N$, $U\ge N$ are natural numbers, $x>x_0$, $y$ are real numbers, $c$ is an absolute constant, not always the same.

Then Haselgrove~S.B.~\cite{Haselgrove}, Statulyavichus~V.~\cite{Statiyavichus}, Pan~C.D. and Pan~C.B.~\cite{Pan Cheng-dong-Pan Cheng-biao}, Zhan~T.~\cite{Zhan Tao-1991} for the sum $S_1(\alpha;x,y)$, $y\ge x^{\theta}$, having obtained a nontrivial estimate on the minor arcs and studied its behavior on the major arcs, proved an asymptotic formula in the ternary Goldbach problem with almost equal summands with conditions $|p_i-N/3|\le H$, $H=N^{\theta}$, respectively for
$$
\theta=\frac{63}{64}+\varepsilon, \qquad  \frac{279}{308}+\varepsilon, \qquad  \frac{2}{3}+\varepsilon ,\qquad  \frac{5}{8}+\varepsilon.
$$

Liu~J. and Zhan~T.~\cite{Liu+Zhan-Mh.Math1999}, having obtained a nontrivial estimate of the sum $S_2(\alpha;x,y)$ both on the minor and on the major arcs for $y\ge x^{\frac{11}{16}+\varepsilon}$, proved the theorem that a sufficiently large natural number $N$ can be represented in the form
$$
N = p_1+p_2+ p_3^2, \qquad \left|p_j-\frac{N}{3}\right|\le H, \qquad \left|p_3^2-\frac N3\right|\le H, \qquad H\ge N^{\frac{27}{32}+\varepsilon},
$$
and in the papers \cite{Liu+Zhan-AcMathSin2000,Liu+Lu+Zhan-Sci-in-China2006} they solved Hua's problem \cite{Hua1938} on the representability of a sufficiently large natural number as a sum of five squares of almost equal primes and showed that a sufficiently large natural number $N$, $N\equiv 5\pmod{24}$, can be represented in the form
$$
N = p_1^2+ \ldots + p_5^2, \quad \left|p_j-\sqrt{\frac{N}{5}}\right|\le H, \quad H\ge N^{\frac{9}{20}+\varepsilon}.
$$

Liu~J. and Zhan~T.~\cite{Liu+Zhan-Konf-Halberstam1996}, using M. Jutila's theorem \cite{Jitila-1991} on the fourth moment of Dirichlet $L$-functions on the critical line, obtained a nontrivial estimate of the sum $S_k(\alpha;x,y)$, $k\ge 3$, on the major arcs $\mathfrak{M}(\mathscr{L}^{c_1})$, $\tau=y^{2k-1}x^{-k+1}\mathscr{L}^{-c_3}$ for $y\ge x^{1-\frac1{2k-1}+\varepsilon}$, where $\mathscr{L}=\ln xq$. Kumchev~A.V. \cite{Kumchev2012} obtained a nontrivial estimate of the sum $S_k(\alpha ;x,y)$ on the minor arcs $\mathfrak{m}(P)$, $\tau=x^{k-\frac{2}{2k+3}}P^{-1}$ for $y\ge x^{1-\frac{1}{2k+3}+\varepsilon}$.
Yao~Y.~\cite{Yanjun2014}, using these estimates, generalized Hua's theorem \cite{Hua1938} in the Waring--Goldbach problem for cubes, i.e., proved that every sufficiently large odd natural number $N$ can be represented in the form
$$
p_1^3+p_2^3+\ldots+p_9^3=N, \qquad \left|p_i-\sqrt[3]{\frac{N}{9}}\right|\le N^{\frac{1}{3}-\frac{1}{51}+\varepsilon}.
$$

In 2016, Z.Kh. Rakhmonov and F.Z. Rakhmonov \cite{RakhmonovZKh+RFZ-DANRT-59-2016,RakhmonovZKh+RFZ-TrMIAN-296-2016}, using I.M. Vinogradov's method of estimating trigonometric sums with primes, obtained a nontrivial estimate of the form
$$
S_3(\alpha ;x,y)\ll \frac{y}{\mathscr{L}^B},
$$
on the minor arcs $\mathfrak{m}\left(\mathscr{L}^{32(B+20)}\right)$, $\tau=y^5x^{-2}\mathscr{L}^{-32(B+20)}$ for $y\ge x^\frac{4}{5}\mathscr{L}^{8B+151}$, where $B$ is an absolute constant (see also \cite{Rakhmonovi-DoklANRT2013,Rakhmonovi-DoklMath2014_124,RakhmonovZamonov2014,RZH+RFZ+ZBM-ChSbor-2016}.

The main result of this paper is a theorem on the behavior of the sum $S_k(\alpha;x,y)$ on the major arcs $\mathfrak{M}(\mathscr{L}^b)$, $\tau=y^{2k-1}x^{-k+1}\mathscr{L}^{-b_1}$, where $b$, $b_1$ are arbitrary fixed positive numbers.

{\theorem\label{Teorema povedenie kor trigsumm s prost chislami} Let $x\ge x_0$, $k\ge3$ be a fixed integer, $A$, $b_1$, $b$ be arbitrary fixed positive numbers, $1\le q\le\mathscr{L}^b$,
$$
\alpha=\frac{a}{q}+\lambda,\quad (a,q)=1,\quad |\lambda|\le\frac{1}{q\tau},\quad \tau=\frac{y^{2k-1}}{x^{k-1}}\mathscr{L}^{-b_1}.
$$
Then for $\lambda\le\left(2\pi k^2x^{k-2}y^2\right)^{-1}$ and $y\ge x^{\frac{\scriptstyle5}{\scriptstyle8}}\mathscr{L}^{1.5A+0.25b+18}$ the equality
$$
S_k(\alpha;x,y)=\frac{1}{\varphi(q)}\sum_{\substack{n=1\\(a,q)=1}}^qe\left(\frac{an^k}{q}\right)\int_{x-y}^xe(\lambda u^k)du+O\left(y\mathscr{L}^{-A}\right),
$$
holds, and for $\lambda>\left(2\pi k^2x^{k-2}y^2\right)^{-1}$ and $y\ge x^{1-\frac{1}{2k-1+\eta_k}}\mathscr{L}^{c_k}$, where $$
\eta_k=\frac{2}{4k-5+2\sqrt{(2k-2)(2k-3)}}, \qquad c_k=
 \frac{2A+22+\left(\frac{2\sqrt{2k-3}}{\sqrt{2k-2}}-1\right)b_1}{2\sqrt{(2k-2)(2k-3)}-(2k-3)},
$$
the estimate
$$
S_k(\alpha;x,y)\ll y\mathscr{L}^{-A}.
$$ holds. }

The explicit values of the parameter $\eta_k$ for $k$ equal to $3$, $4$, $5$ and $6$ have the form:
\begin{align*}
&\eta_3=\frac2{7+4\sqrt{3}}\approx0.1435935394,
&\eta_4=\frac2{11+2\sqrt{30}}\approx0.0910976998,\\
&\eta_5=\frac2{15+4\sqrt{14}}\approx0.0667409058,
&\eta_6=\frac2{19+2\sqrt{90}}\approx0.052668078.
\end{align*}

This theorem is slightly weaker than the above-mentioned theorem of Liu~J. and Zhan~T.~\cite{Liu+Zhan-Konf-Halberstam1996}, and its proof uses the theorem on the second moment of Dirichlet $L$-functions on the critical line \cite{RanVV-J.Lon.Math.Soc-1980,Zhan Tao-Acta-Math-Sinica-1992} (see also \cite{RakhmonovZKh+SAA-FPM-DANRT-63(5-6)-2020, RakhmonovZKh+SAA-FPM-DANRT-63(7-8)-2020}).

\section{Known Lemmas}
 {\lemma\label{Lemma otsenka trig integ po pervoy proizv} Let the real function $f(u)$ and the monotone function $g(u)$ satisfy the conditions: $f'(u)$ is monotone, $|f'(u)|\ge m_1>0$ and $|g(u)|\le M$. Then the estimate holds:
$$
\int_a^bg(u)e(f(u))du\ll \frac{M}{m_1}.
$$}

{\sc Proof} see \cite{AGI+KAA+ChVN-TKTS-1987-Nauka}

 {\lemma\label{Lemma otsenka trig integ po vtoroy proizv} Let the real function $f(u)$ and the monotone function $g(u)$ satisfy the conditions: $f'(u)$ is monotone, $|f''(u)|\ge m_2>0$ and $|g(u)|\le M$. Then the estimate holds:
$$
\int_a^bg(u)e(f(u))du\ll \frac{M}{\sqrt{m_2}}.
$$}

{\sc Proof} see \cite{AGI+KAA+ChVN-TKTS-1987-Nauka}

{\lemma\label{Lemma razlozh psi po nuyam L-ryadov} Let $2\le T_0\le x$, $\rho =\beta +i\gamma$ be nontrivial zeros of the function $L(s,\chi)$. Then
$$
\psi(x,\chi)=E_0x-\sum_{|\gamma |\le T_0}\frac{x^{\rho}}{\rho}+R(x,T_0),\quad R(x,T_0)\ll \frac{x\mathscr{L}^2}{T_0},
$$
where $E_0=1$ if $\chi =\chi_0$; $E_0=0$ if $\chi \ne \chi _0$.}

{\sc Proof} see \cite{Devenport-MTCh-1971-Nauka}

 {\lemma\label{Lemma granitsa nuley L-ryada} For a suitable $c>0$ the function $L(s,\chi )$, $s=\sigma+it$, has no zeros in the region
$$
\sigma\ge 1-\delta(q,t),\quad \delta(q,t)=\frac{c}{\max(\ln q,\ln^{3/4}(|t|+3)\ln^{3/4}\ln (|t|+3))},
$$
for all characters $\chi$ modulo $q$, except possibly for a simple real zero $\beta_1$ of the $L$-function defined by an exceptional character $\chi _1$.}

{\sc Proof} see \cite{Prachar-RPCh-1967-Mir}

{\lemma \label{Lemma kol-va nuley L-ryada v pryamougolnike} The number of zeros $\rho$ of the function $L(s,\chi )$, $\chi\pmod q$, for which $T\le|\gamma|\le T+1$, does not exceed $c\ln qT$.}

{\sc Proof} see \cite{Kasratsuba-OATCh-1983-Nauka}.

{\lemma\label{Lemma granitsa deystvitel nulya L-ryada} For any $\varepsilon >0$ there exists $c=c(\varepsilon )$ such that if $\chi_1$ is a real character modulo $q$ and $\beta_1$ is a real zero of $L(s,\chi_1)$, then
$$
\beta_1\le 1-\frac {c(\varepsilon)}{q^{\varepsilon}}.
$$}

{\sc Proof} see \cite{Kasratsuba-OATCh-1983-Nauka}.

{\lemma\label{Lemma plotnostn-teor Zhan Tao} Let $\varepsilon$ be an arbitrarily small positive constant, and $T^{\frac{35}{108}+\varepsilon}\le H\le T$, then the estimates hold
\begin{align*}
&\sum_{\chi\bmod q}\left(N(u,T+H,\chi)-N(u,T,\chi)\right)\ll\left\{
                                                       \begin{array}{ll}
                                                         (qH)^{\frac{4}{3-2u}(1-u)}(\ln qH)^9, & \hbox{for} \quad\dfrac{1}{2}\le u\le\dfrac{3}{4} \vspace{7pt} \\
                                                         (qH)^{\frac{2}{u}(1-u)+\varepsilon}, & \hbox{for}\quad\dfrac{3}{4}\le u\le 1,
                                                       \end{array}
                                                     \right.
 \end{align*}}

{\sc Proof} see \cite{Zhan Tao-Acta-Math-Sinica-1992}.

{\lemma \label{Lemma otseka smeshann summi} Let $a$, $k$ and $q$ be natural numbers, $(a,q)=1$, $\nu_k(q)$ be the number of $k$-th power residues modulo $q$,
$$
\tau(\chi,a,k)=\sum_{h=1}^q\chi(h)e\left(\frac{ah^k}q\right),
$$
Then the inequality holds
$$
\left|\tau(\chi,a,k)\right|\le 2(\tau(q))^{\frac{\ln k}{\ln2}}\sqrt{q}.
$$}

{\sc Proof} The assertion of the lemma follows from problems $8.\alpha)$ and $8.\beta)$, p. 103 of the textbook \cite{Vinogradov Osn-Teor-Chisel-1981}.

Here is the English translation of the LaTeX mathematical text:

\section{Proof of Theorem \ref{Teorema povedenie kor trigsumm s prost chislami}}
We divide the proof of the theorem into several stages.

\subsection{Reduction of the proof of the theorem to the estimate of the sum $W(\alpha;x,y)$}
From the condition $1\le q\le\mathscr{L}^b$, $k\ge3$ a fixed integer, it follows that
$$
\frac{q\left(\tau(q)\right)^\frac{\ln k}{\ln 2}}{\varphi (q)}\ll\mathscr{L}.
$$
Applying Lemma \ref{Lemma otseka smeshann summi} to the sum $\tau(\chi,a,k)$, and then the last estimate, we find
\begin{equation}\label{Formula otsenka tau(chi,a,k)}
|\tau(\chi,a,k)|=\left|\sum_{h=1}^q\chi(h)e\left(\frac{ah^k}q\right)\right|\le 2\sqrt{q}\left(\tau(q)\right)^\frac{\ln k}{\ln 2}\ll\frac{\varphi(q)\mathscr{L}}{\sqrt{q}}.
\end{equation}
Using the orthogonality property of characters, we find
\begin{align}
S_k(\alpha;x,y)&=\frac{1}{\varphi(q)}\sum_{\chi\hspace{-7pt}\mod q}\tau(\bar{\chi},a,k)\sum_{x-y<n\le x}\Lambda(n) \chi(n)e(\lambda n^k)+O(\mathscr{L}^2). \label{formula Sk(alpha ;x,y)=summchi...}
\end{align}
Further, applying Abel's transformation in integral form, we have
\begin{eqnarray*}
\sum_{x-y<n\le x}\Lambda(n)\chi(n)e(\lambda n^k)=-\int_{x-y}^{x}\psi(u,\chi)de(\lambda u^k)+e(\lambda x^k)\psi(x,\chi)-e(\lambda(x-y)^k)\psi (x-y,\chi).
\end{eqnarray*}
Using Lemma \ref{Lemma razlozh psi po nuyam L-ryadov} with $T_0=\left(xy^{-1}+|\lambda|x^k\right)q^{\frac12}\mathscr{L}^{A+3}$, we find
\begin{align*}
\sum_{x-y<n\le x}\Lambda(n)\chi(n)e(\lambda n^k)=-&\int_{x-y}^x\left(E_0u- \sum_{|\gamma |\le T_0}\frac{u^\rho}{\rho}\right)de(\lambda u^k)+
e(\lambda x^k)\left(E_0x-\sum_{|\gamma|\le T_0}\frac{x^\rho}{\rho}\right)-\\
&-e(\lambda(x-y)^k)\left(E_0(x-y)-\sum_{|\gamma|\le T_0}\frac{(x-y)^\rho}{\rho}\right)-\\
-\int_{x-y}^xR(u,T_0)2\pi ik\lambda&u^{k-1}e(\lambda u^k)du+e(\lambda x^k)R(x,T_0)-e(\lambda(x-y)^k)R(x-y,T_0).
\end{align*}
Applying integration by parts to the first integral, and also using the estimate for $R(u;T_0,\chi )$ from Lemma \ref{Lemma razlozh psi po nuyam L-ryadov} and the value of the parameter $T_0$, we find
\begin{align*}
\sum_{x-y<n\le x}\Lambda(n)\chi(n)e(\lambda n^k)&=E_0\int_{x-y}^xe(\lambda u^k)du-\sum_{|\gamma|\le T_0} I(\rho,\lambda)+ O\left(\frac{y}{q^\frac12\mathscr{L}^{A+1}}\right),\\
&I(\rho,\lambda)
=\int_{x-y}^xu^{\beta-1}e\left(\lambda u^k+\frac\gamma{2\pi}\ln u\right)du.
\end{align*}
Substituting the found formula into (\ref{formula Sk(alpha ;x,y)=summchi...}) and using estimate (\ref{Formula otsenka tau(chi,a,k)}), we find
\begin{align}
S_k(\alpha ;x,y)&=\frac{\tau(\chi_0,a,k)}{\varphi(q)}\int_{x-y}^xe(\lambda u^k)du-W(\alpha;x,y)-E_1 W_1(\alpha;x,y)+O(y\mathscr{L}^{-A}), \label{formula Sk(alpha;x,y)=int+W1+W+R2} \\
&W(\alpha;x,y)=\frac{1}{\varphi(q)}\sum_{\chi\hspace{-7pt}\mod q}\tau(\bar{\chi},a,k)\sum_{\substack{|\gamma|\le T_0\\\rho\ne\beta_1}}I(\rho,\lambda),\nonumber \\
&W_1(\alpha;x,y)=\frac{\tau(\chi_1,a,k)}{\varphi(q)}I(\beta_1,\lambda), \nonumber
\end{align}
where $E_{1}=1$ if modulo $q$ there exists a real character $\chi_1$ such that $L(s,\chi_1)$ has a real zero $\beta_1$, $\beta_1\ge 1-c/\ln q$, and $E_1=0$ otherwise.

\textbf{Estimate of $|W_1(\alpha;x,y)|$}. Using trivial estimates of the sum $\tau(\chi_1,a,k)$ and the integral $I(\beta_1,\lambda)$, we find
$$
|W_1(\alpha;x,y)|=\left|\frac{\tau(\chi_1,a,k)}{\varphi(q)}\int_{x-y}^x\frac{u^{\beta_1-1}}{\beta_1}e(\lambda u^k)du\right|\ll yx^{\beta_1-1}.
$$
According to Lemma \ref{Lemma granitsa deystvitel nulya L-ryada}, bearing in mind that $q\le\mathscr{L}^b$, for $\varepsilon=(2b)^{-1}$ we have
\begin{align*}
x^{\beta_1-1}=\exp\left((\beta_1-1)\mathscr{L})\right)\le\exp\left(-\frac{c(\varepsilon)\mathscr{L}}{q^\varepsilon}\right)\le \exp\left(-\frac{c(\varepsilon )\mathscr{L}}{\mathscr{L}^{b\varepsilon}}\right)=\exp(-c(\varepsilon)\sqrt{\mathscr{L}}).
\end{align*}
Consequently,
\begin{equation}\label{Formula otsenka summi W1}
|W_1(\alpha;x,y)|\ll y\exp(-c(\varepsilon)\sqrt{\mathscr{L}})\ll y\mathscr{L}^{-A}.
\end{equation}

\textbf{Transformation of $|W(\alpha;x,y)|$}. Passing to estimates, we have
\begin{equation}\label{Formula W(alpha,x,y)<=frac(1)(varphi(q))summ|tau()|.W(lambda)}
|W(\alpha;x,y)|\le\frac{1}{\varphi(q)}\sum_{\chi\hspace{-7pt}\mod q}|\tau(\bar{\chi},a,k)|\EuScript{W}(\lambda,\chi),\qquad
\EuScript{W}(\lambda,\chi)=\sum_{\substack{|\gamma|\le T_0\\\rho\ne\beta_1}}\left|I(\rho,\lambda)\right|,
\end{equation}
where $\beta_1$ is a real zero if modulo $q$ there exists a real character $\chi_1$ such that $L(s,\chi_1)$ has a real zero $\beta_1\ge 1-c/\ln q$. We shall estimate the sum $|W(\alpha;x,y)|$ only in the case $\lambda\ge 0$. The case $\lambda\le 0$ reduces to the case $\lambda\ge 0$ by means of the relation
\begin{align*}
\EuScript{W}(\lambda,\chi)&=\sum_{\substack{|\gamma|\le T_0\\\rho\ne\beta_1}}\left|\overline{I(\rho,\lambda)}\right|=\sum_{\substack{|\gamma|\le T_0\\\rho\ne\beta_1}}\left|\int_{x-y}^{x}u^{\beta -1}e\left(-\lambda u^k-\frac{1}{2\pi}\gamma \ln u\right)du\right|=\\
&=\sum_{\substack{|\gamma|\le T_0\\\rho\ne\beta_1}}\left|I(\bar{\rho},-\lambda)\right|=\sum_{\substack{|\gamma|\le T_0\\\rho\ne\beta_1}}\left|I(\rho,-\lambda)\right|=\EuScript{W}(\chi,-\lambda),
\end{align*}
Estimating the integral $I(\rho,\lambda)$ using Lemma \ref{Lemma otsenka trig integ po pervoy proizv}, with $M=x^{\beta-1}$, $f(u)=\lambda u^k+\frac{\gamma}{2\pi}\ln u$ and $m_1=\min f'(u)$, we find
\begin{align*}
&|I(\rho,\lambda)|\le\frac{x^{\beta-1}}{\min|f'(u)|}=\frac{x^\beta}{\min|xf'(u)|},&f'(u)=k\lambda u^{k-1}+\frac{\gamma}{2\pi u}= \frac{\gamma+2\pi k\lambda u^k}{2\pi u}.
\end{align*}
We also estimate the integral $I(\rho,\lambda)$ using Lemma \ref{Lemma otsenka trig integ po vtoroy proizv}, setting $M=x^{\beta-1}$, $m_2=\min f''(u)$. We have
\begin{align*}
&|I(\rho,\lambda)|\le\frac{x^{\beta-1}}{\sqrt{\min|f''(u)|}}=\frac{x^\beta}{\sqrt{\min|x^2f''(u)|}},&f''(u)=\frac{2\pi k(k-1)\lambda u^k-\gamma}{2\pi u^2}.
\end{align*}
Consequently, taking into account the following trivial estimate
$$
|I(\rho,\lambda)|\le\int_{x-y}^{x}u^{\beta-1}du\le yx^{\beta-1},
$$
we have
$$
|I(\rho,\lambda)|\le x^\beta\min\left(\frac{y}{x},\frac{1}{\min|xf'(u)|},\frac{1}{\sqrt{\min|x^2f''(u)|}}\right).
$$
Using the relations
\begin{align*}
&|xf'(u)|=|\gamma+2\pi k\lambda u^k|\frac{x}{2\pi u}\ge|\gamma+2\pi k\lambda u^k|\frac{x}{2\pi x}\ge\frac1{2\pi}|\gamma+2\pi k\lambda u^k|,\\
&|x^{2}f''(u)|=\frac{|2\pi k(k-1)\lambda u^k-\gamma|x^2}{2\pi u^2}\ge\frac{|2\pi k(k-1)\lambda u^k-\gamma|}{2\pi },
\end{align*}
we represent this estimate in the form
\begin{align}
|I(\rho,\lambda)|\ll x^\beta\min\left(\frac{y}{x},\frac{1}{\min|\gamma+2\pi k\lambda u^k|},\frac{1}{\sqrt{|2\pi k(k-1)\lambda u^k-\gamma|}}\right). \label{Formula otsenka I(rho,lambda)}
\end{align}
Substituting this estimate and (\ref{Formula otsenka tau(chi,a,k)}) into (\ref{Formula W(alpha,x,y)<=frac(1)(varphi(q))summ|tau()|.W(lambda)}), we obtain
\begin{align}\label{Formula W(lambda,x,y)<<..sum W(lambda,chi)}
&|W(\alpha;x,y)|\ll\frac{\mathscr{L}}{\sqrt{q}}\sum_{\chi\hspace{-7pt}\mod q}\EuScript{W}(\lambda,\chi).
\end{align}

Here is the English translation of the LaTeX mathematical text:

\subsection{{Proof of the theorem in the case $\lambda\le\left(2\pi k^2x^{k-2}y^2\right)^{-1}$}}
In this subsection, without loss of generality, we shall assume that
\begin{equation}\label{Formula usl <=y>=sl1}
y=x^{\frac{\scriptstyle5}{\scriptstyle8}}\mathscr{L}^{1.5A+0.25b+18}.
\end{equation}

We divide all zeros $\rho=\beta+i\gamma$ of the function $L(s,\chi)$ with the condition $|\gamma|\le T_0$ into sets $D_1$, $D_2$ and $D_3$ as follows:
\begin{align*}
&D_1=\left\{\rho:\ -T_0\le\gamma<-2\pi k\lambda x^k-\frac{x}{y}\right\},  \\
&D_2=\left\{\rho:\ -2\pi k\lambda x^k-\frac{x}{y}\le\gamma\le-2\pi k\lambda(x-y)^k+\frac{x}{y}\right\}, \\
& D_3=\left\{\rho:\ -2\pi k\lambda(x-y)^k+\frac{x}{y}<\gamma\le T_0 \right\}.
\end{align*}
Denoting by $\EuScript{W}_j(\lambda,\chi)$, $j=1,\ 2,\ 3$, the sum of the moduli of the integral $I(\rho,\lambda)$ over the zeros $\rho$ belonging to the set $D_j$, we represent the sum $\EuScript{W}(\lambda,\chi)$ in (\ref{Formula W(alpha,x,y)<=frac(1)(varphi(q))summ|tau()|.W(lambda)}) in the form:
\begin{align}\label{Formula W(chi,lambda)=W1+W2+W3,lambda<}
&\EuScript{W}(\lambda,\chi)=\EuScript{W}_1(\lambda,\chi)+\EuScript{W}_2(\lambda,\chi)+\EuScript{W}_3(\chi,\lambda).
\end{align}

\textbf{Estimate of $\EuScript{W}_1(\lambda,\chi)$.} Adding the term $2\pi k\lambda u^k$, $x-y\leq u \le x$, to all three members of the inequality defining the set $D_1$, we obtain
\begin{align*}
&D_1=\left\{\rho:\ -T_0+2\pi k\lambda u^k\le\gamma+2\pi k\lambda u^k<-2\pi k\lambda x^k+2\pi k\lambda u^k-\frac{x}{y}\right\},  \\
\end{align*}
On the interval $x-y\leq u\le x$ the function $2\pi k\lambda u^k$ is monotonically increasing, therefore for the right boundary of the set $D_1$ we have
$$
-2\pi k\lambda x^k+2\pi k\lambda u^k-\frac{x}{y}\le-\frac{x}{y}.
$$
Consequently, if $\rho$ belongs to the set $D_1$, then the inequality $\gamma+2\pi k\lambda u^k<-\frac{x}{y}$ holds, therefore for the monotonically increasing function $\gamma+2\pi k\lambda u^k$ on the interval $x-y\le u\le x$ the relation
\begin{align*}
&\min|\gamma+2\pi k\lambda u^k|=-\max(\gamma+2\pi k\lambda u^k)=-\gamma-2\pi k\lambda x^k\ge\frac{x}{y},
\quad \text{if}  \quad \rho\in D_1.
\end{align*}
holds. Hence, taking into account the second estimate (\ref{Formula otsenka I(rho,lambda)}), we find
$$
\EuScript{W}_1(\lambda,\chi)\ll\sum_{\rho \in D_1}\frac{x^\beta }{-\gamma-2\pi k\lambda x^k}.
$$
We divide all zeros in the set
$$
D_1=\left\{\rho:\ \frac{x}{y}\le-\gamma-2\pi k\lambda x^k<T_0-2\pi k\lambda x^k\right\},
$$
into classes $ D_{11},\ \ldots,\ D_{1r}$, $r\ll T_0x^{-1}y$, as follows: in the class $D_{1n}$ we place those zeros $\rho$ for which the conditions
$$
\frac{nx}y<-\gamma-2\pi k\lambda x^k\le\frac{(n+1)x}y
$$
hold. Therefore
\begin{align*}
\EuScript{W}_1(\lambda,\chi)\ll&\sum_{n=1}^r\sum_{\rho\in D_{1n}}\frac{x^\beta}{-\gamma-2\pi k\lambda x^k} \le\frac{y}{x}\sum_{n=1}^r\sum_{\rho\in D_{1n}}\frac{x^\beta }{n}\le\frac{y\mathscr{L}}x\max_{1\le n\le r}\sum_{\rho\in D_{1n}}\hspace{-6pt}x^{\beta}\le\frac{y\mathscr{L}}x\max_{|T|\le T_0}\sum_{T-\frac{x}y<\gamma\le T}\hspace{-10pt}x^\beta.
\end{align*}

\textbf{Estimate of $\EuScript{W}_3(\chi,\lambda)$.} Adding the term $2\pi k\lambda u^k$, $x-y\leq u\le x$, to all three members of the inequality defining the set $D_3$, we obtain
\begin{align*}
& D_3=\left\{\rho:\ 2\pi k\lambda u^k-2\pi k\lambda(x-y)^k+\frac{x}{y}<\gamma+2\pi k\lambda u^k\le T_0+2\pi k\lambda u^k \right\}.
\end{align*}
On the interval $x-y\leq u\le x$ the function $2\pi k\lambda u^k$ is monotonically increasing, therefore for the left boundary of the set $D_3$ we have
$$
2\pi k\lambda u^k-2\pi k\lambda(x-y)^k+\frac{x}{y}\ge\frac{x}{y}.
$$
Consequently, if $\rho$ belongs to the set $D_3$, then the inequality $\gamma+2\pi k\lambda u^k>\frac{x}{y}$ holds, therefore for the monotonically increasing function $\gamma+2\pi k\lambda u^k$ on the interval $x-y\le u\le x$ the relation
\begin{align*}
&\min|\gamma+2\pi k\lambda u^k|=\min(\gamma+2\pi k\lambda u^k)=\gamma+2\pi k\lambda(x-y)^k\ge\frac{x}{y},
\quad \text{if}\quad\rho\in D_3.
\end{align*}
holds. Hence, taking into account the second estimate (\ref{Formula otsenka I(rho,lambda)}), we find
$$
\EuScript{W}_3(\chi,\lambda)\ll\sum_{\rho \in D_3}\frac{x^\beta }{\gamma+2\pi k\lambda(x-y)^k}.
$$
We divide all zeros in the set
$$
 D_3=\left\{\rho:\ \frac{x}{y}\le\gamma+2\pi k\lambda(x-y)^k<T_0+2\pi k\lambda(x-y)^k\right\},
$$
into classes $D_{31},\ \ldots,\ D_{3r}$, $r\ll T_0x^{-1}y$, as follows: in the class $D_{3n}$ we place those zeros $\rho$ for which the conditions
$$
\frac{nx}{y}<\gamma+2\pi k\lambda(x-y)^k\le\frac{(n+1)x}{y},\qquad\text{if}\qquad 1\le n\le r .
$$
hold. Therefore
\begin{align*}
\EuScript{W}_3(\chi,\lambda)&\ll\sum_{n=1}^r\sum_{\rho\in D_{3n}}\frac{x^\beta}{\gamma+2\pi k\lambda(x-y)^k}\le\frac{y}{x}\sum_{n=1}^r\sum_{\rho\in D_{3n}}\frac{x^\beta }{n}\le\frac{y\mathscr{L}}x\max_{|T|\le T_0}\sum_{T-\frac{x}y<\gamma\le T}x^\beta.
\end{align*}

\textbf{Estimate of $\EuScript{W}_2(\lambda,\chi)$.} Representing the set $D_2$ in the form
$$
D_2=\left\{\rho:\ T_1-2\pi k\lambda(x^k-(x-y)^k)+\frac{2x}{y}\le-\gamma\le T_1\right\}, \qquad T_1=2\pi k\lambda x^k+\frac{x}{y}\le T_0,
$$
and bearing in mind that for $\lambda\le\left(2\pi k^2x^{k-2}y^2\right)^{-1}$ for the length of the set $D_2$ the inequality
$$
2\pi k\lambda(x^k-(x-y)^k)+\frac{2x}{y}\le2\pi k^2\lambda x^{k-1}y+\frac{2x}{y}\le\frac{(2k+2)x}{y},
$$
holds, and using the trivial estimate of the integral $I(\rho,\lambda)$, i.e. the first estimate (\ref{Formula otsenka I(rho,lambda)}), we have
$$
\EuScript{W}_2(\lambda,\chi)\le\sum_{\rho \in D_2}|I(\rho,\lambda)|\ll\frac{y}{x}\sum_{\rho \in D_2}x^\beta\le\frac{(2k+2)y}{x}\max_{|T|\le T_0}\sum_{T-\frac{x}{y}\le-\gamma\le T}x^\beta\ll\frac{y}{x}\max_{|T|\le T_0}\sum_{T-\frac{x}{y}\le\gamma\le T}x^\beta.
$$
Substituting the obtained estimates for $\EuScript{W}_1(\lambda,\chi)$, $\EuScript{W}_3(\chi,\lambda)$ and $\EuScript{W}_3(\chi,\lambda)$ into (\ref{Formula W(chi,lambda)=W1+W2+W3,lambda<}), and then into (\ref{Formula W(lambda,x,y)<<..sum W(lambda,chi)}),
we find
\begin{align}
&|W(\alpha;x,y)|\ll\frac{y\mathscr{L}^2}{\sqrt{q}\,x}\max_{|T|\le T_0}\EuScript{V}_q\left(T,\frac{x}y\right), \label{Formula W(alpha,x,y)<<frac()()...W(lambda)}\\
&\EuScript{V}_q(T,U)=\sum_{\chi\hspace{-7pt}\mod q}\ \sum_{T-U<\gamma\le T}x^\beta, \qquad U\le T. \nonumber
\end{align}
To estimate $\EuScript{V}_q(T,U)$ we use the density theorem in narrow rectangles of the critical strip for zeros of Dirichlet $L$-series modulo $q$ and the theorem on the boundary of these zeros. We have
\begin{align*}
\EuScript{V}_q(T,U)=&\mathscr{L}\int_0^1x^u\sum_{\chi\hspace{-7pt}\mod q}\left(N(u,T,\chi)-N(u,T-U,\chi)\right)du+\sum_{\chi\hspace{-7pt}\mod q}\left(N(T,\chi)-N(T-H,\chi)\right)\le\\
&\ll\mathscr{L}\max_{0.5\le u\le1}x^u\sum_{\chi\hspace{-7pt}\mod q}\left(N(u,T,\chi)-N(u,T-U,\chi)\right).
\end{align*}
According to Lemma \ref{Lemma granitsa nuley L-ryada} the function $L(u+it,\chi )$ has no zeros in the region
$$
u\ge1-\delta(q,t),\qquad \delta(q,t)=\frac{c_1}{\max\left(\ln q,\ \left(\ln(t+3)\ln\ln(t+3)\right)^\frac{3}{4}\right)},
$$
for all characters $\chi\mod q$, except possibly for a simple real zero $\beta_1$ of the $L$-function defined by an exceptional character $\chi _1$. Therefore, bearing in mind that $\delta(q,T)\ge\delta(q,T_0)$, we find
\begin{equation}\label{Formula otsenka Vq(T,U)}
\EuScript{V}_q\left(T,U\right)\ll\mathscr{L}\max_{0.5\le u\le1-\delta}x^u\sum_{\chi\hspace{-7pt}\mod q}\left(N(u,T,\chi)-N(u,T-U,\chi)\right), \qquad \delta=\delta(q,T_0).
\end{equation}
Substituting the right-hand side of this inequality for $U=\frac{x}{y}$ into (\ref{Formula W(alpha,x,y)<<frac()()...W(lambda)}), we find
\begin{align}
|W(\alpha;x,y)|\ll\frac{y\mathscr{L}^3}{\sqrt{q}\,x}\max_{|T|\le T_0}\max_{0.5\le u\le1-\delta}x^u\sum_{\chi\hspace{-7pt}\mod q}\left(N(u,T,\chi)-N\left(u,T-\frac{x}{y},\chi\right)\right).
\label{Formula W(lambda,x,y)<<..maxmax..sum(N(u,T,chi)-N(u,T-fracxy,chi))}
\end{align}
From the relation $|T|\le T_0=\left(xy^{-1}+|\lambda|x^k\right)q^{\frac12}\mathscr{L}^{A+3}$, $0\le\lambda\le\left(2\pi k^2x^{k-2}y^2\right)^{-1}$ and condition (\ref{Formula usl <=y>=sl1}), we have
\begin{align*}
\frac{T}{(xy^{-1})^3}\le&\left(\frac{y^2}{x^2}+\lambda x^{k-3}y^3\right)q^{\frac12}\mathscr{L}^{A+3}\le\left(\frac{y^2}{x^2}+\frac{y}{2\pi k^2x}\right)q^{\frac12}\mathscr{L}^{A+3}\le1.
\end{align*}
Hence it follows that in the sum over $\chi\hspace{-5pt}\mod q$ in (\ref{Formula W(lambda,x,y)<<..maxmax..sum(N(u,T,chi)-N(u,T-fracxy,chi))}) the condition $\frac{x}{y}\ge T^\frac13$ holds, i.e. Lemma \ref{Lemma plotnostn-teor Zhan Tao} can be applied to this sum. Dividing the right-hand side of (\ref{Formula W(lambda,x,y)<<..maxmax..sum(N(u,T,chi)-N(u,T-fracxy,chi))}) the interval $0.5\le u\le1-\delta$ into the intervals $0.5\le u\le 0.75$ and $0.75\le u\le1-\delta$ we obtain two sums, which we shall denote by $\mathscr{A}_1$ and $\mathscr{A}_2$, respectively, and then setting in Lemma \ref{Lemma plotnostn-teor Zhan Tao} $\varepsilon=\frac{\delta}{3}$, $\delta=\delta(q,T_0)$, we have
\begin{align}
&|W(\alpha;x,y)|\ll\mathscr{A}_1+\mathscr{A}_2,\label{Formula W(alpha;x,y)<<A1+A2}\\
&\mathscr{A}_1=\frac{y\mathscr{L}^{12}}{\sqrt{q}\,x}\max_{0.5\le u\le0.75}x^u\left(\frac{qx}{y}\right)^{\frac{4-4u}{3-2u}},\nonumber \\
&\mathscr{A}_2=\frac{y\mathscr{L}^3}{\sqrt{q}\,x}\max_{0.75\le u\le1-\delta}x^u\left(\frac{qx}{y}\right)^{\frac{2}{u}(1-u)+\varepsilon}. \nonumber
\end{align}
Note that the values of $\mathscr{A}_1$ and $\mathscr{A}_2$ at $u=0.75$ coincide by definition, therefore, for uniformity, after applying Lemma \ref{Lemma plotnostn-teor Zhan Tao}, for the upper estimate of $\mathscr{A}_1$ and $\mathscr{A}_2$ at $u=0.75$ we take the smaller quantity, i.e. the estimate for $\mathscr{A}_1$.

\textbf{Estimate of $\mathscr{A}_1$}. We have
\begin{align*}
&\mathscr{A}_1=\frac{q^{1.5}x\mathscr{L}^{12}}{y}\max_{0.5\le u\le0.75}f_1(u),\qquad f_1(u)=x^u\left(\frac{qx}{y}\right)^\frac1{u-1.5}>0,\\
&f_1'(u)=f_1(u)\left(\ln x+\frac{\ln\left(\frac{y}{qx}\right)}{(u-1.5)^2}\right)=\frac{f_1(u)}{(u-1.5)^2}\ln\frac{y}{qx^{-u^2+3u-1.25}}.
\end{align*}
From conditions (\ref{Formula usl <=y>=sl1}) and $q\le\mathscr{L}^b$, as well as from the relation
$$
\max_{0.5\le u\le0.75}\left(-u^2+3u-1.25\right)=\left.\left(-u^2+3u-1.25\right)\right|_{u=0.75}=\frac{7}{16},
$$
it follows that
$$
\ln\frac{y}{qx^{-u^2+3u-1.25}}\ge\ln\frac{yx^\frac58\mathscr{L}^{1.5A+0.25b+18}}{\mathscr{L}^bx^\frac7{16}}\ge\ln x^\frac3{16}\mathscr{L}^{1.5A+18-0.75b}>\ln x^\frac18>0,
$$
i.e. $f_1'(u)>0$ and $f_1(u)$ is an increasing function. Using this property, and then relation (\ref{Formula usl <=y>=sl1}), we have
\begin{align*}
&\mathscr{A}_1=\frac{q^{1.5}x\mathscr{L}^{12}}{y}f_1(0.75)=
y\left(\frac{x^\frac58q^\frac14\mathscr{L}^{18}}{y}\right)^\frac23\le y\left(\frac{x^\frac58\mathscr{L}^{0.25b+18}}{y}\right)^\frac23 \ll y\mathscr{L}^{-A}.
\end{align*}

\textbf{Estimate of $\mathscr{A}_2$}. We have
\begin{align*}
\mathscr{A}_2=\frac{y^{3-\varepsilon}\mathscr{L}^3}{q^{2.5-\varepsilon}x^{3-\varepsilon}}\max_{0.75\le u\le1-\delta}f_2(u),\qquad f_2(u)=x^u\left(\frac{qx}{y}\right)^\frac2u>0,\qquad \delta=\delta(q,T_0).
\end{align*}
Hence, using condition (\ref{Formula usl <=y>=sl1}), we have $y=x^{\frac{\scriptstyle5}{\scriptstyle8}}\mathscr{L}^{1.5A+0.25b+18}.$
\begin{align*}
&f_2'(u)=f_2(u)\left(\ln x+\frac{2\ln\left(\frac{y}{qx}\right)}{u^2}\right)=\frac{2f_2(u)}{u^2}\ln\frac{y}{qx^{1-0.5u^2}},\\
&f_2'(0.75)=\frac{32f_2(0.75)}9\ln\frac{y}{qx^\frac{23}{32}}=  \frac{32f_2(0.75)}9\ln\frac{\mathscr{L}^{1.5A+0.25b+18}}{qx^\frac5{32}}<0,\\
&f_2'(1-\delta)=\frac{2f_2(1-\delta)}{(1-\delta)^2}\ln\frac{y}{qx^{0.5+\delta-0.5\delta^2}}=  \frac{2f_2(1-\delta)}{(1-\delta)^2}\ln\frac{x^{\frac18-\delta+0.5\delta^2}\mathscr{L}^{1.5A+0.25b+18}}{q}>0,\\
&f_2''(u)=f_2(u)\left(\ln x+\frac{2\ln\left(\frac{y}{qx}\right)}{u^2}\right)^2-f_2(u)\frac{4\ln\left(\frac{y}{qx}\right)}{u^3}>0,
\end{align*}
Consequently, in the interval $0.75\le u\le1-\delta$ there exists a point $u_0$ such that $f_2'(u_0)=0$ and $ f_2''(u_0)>0$, i.e. in the interval $0.75\le u\le1-\delta$ the graph of the function $f_2(u)$ is concave downward and the sum $\mathscr{A}_2$ is estimated from above by its value at the points $u=0.75$ and $u=1-\delta$. Above we agreed that for the upper estimate of the sum $\mathscr{A}_2$ at $u=0.75$ we take the estimate for $\mathscr{A}_1$ at this point, therefore we have
\begin{align*}
\mathscr{A}_2&\ll y\mathscr{L}^{-A}+\frac{y^{3-\varepsilon}\mathscr{L}^3}{q^{2.5-\varepsilon}x^{3-\varepsilon}}f_2(1-\delta)= y\mathscr{L}^{-A}+\frac{y^{3-\varepsilon}\mathscr{L}^3}{q^{2.5-\varepsilon}x^{3-\varepsilon}}x^{1-\delta}\left(\frac{qx}{y}\right)^\frac2{1-\delta}=\\
&=y\mathscr{L}^{-A}+y\frac{x^{\frac{2\delta}{1-\delta}-\delta+\varepsilon}\mathscr{L}^3}{q^{0.5+\frac{2\delta}{1-\delta}-\varepsilon} y^{\frac{2\delta}{1-\delta}+\varepsilon}}\le y\mathscr{L}^{-A}+y\frac{x^{\frac{2\delta}{1-\delta}-\delta+\varepsilon}\mathscr{L}^3}{y^{\frac{2\delta}{1-\delta}+\varepsilon}}=\\
&=y\mathscr{L}^{-A}+y\mathscr{L}^3\left(\frac{x^{\frac{\scriptstyle5}{\scriptstyle8}}\mathscr{L}^{1.5A+0.25b+18}}y x^{f(\delta,\varepsilon)}\mathscr{L}^{-1.5A-0.25b-18}\right)^{\frac{\scriptstyle{2\delta}}{\scriptstyle{1-\delta}}+\varepsilon},\\
&\qquad f(\delta,\varepsilon)=\frac{\delta+\delta^2+(1-\delta)\varepsilon}{2\delta+(1-\delta)\varepsilon}-\frac58.
\end{align*}
Hence and from relation~\eqref{Formula usl <=y>=sl1}, we obtain
\begin{align*}
\mathscr{A}_2&\ll y\mathscr{L}^{-A}+yx^{f(\delta,\varepsilon)\frac{\scriptstyle{2\delta+(1-\delta)\varepsilon}}{\scriptstyle{1-\delta}}}\mathscr{L}^3.
\end{align*}
Further, for $\varepsilon=\frac{\delta}3$, using the relation
\begin{align*}
f(\delta,\varepsilon)&\frac{2\delta+(1-\delta)\varepsilon}{1-\delta}=
-\frac{\delta}{8}+\frac38\left(\varepsilon-\frac{\delta}{3}-\frac{\delta^2}{3(1-\delta)}\right)\le-\frac{\delta}{8},
\end{align*}
we find
\begin{align*}
\mathscr{A}_2&\ll y\mathscr{L}^{-A}+yx^{-0.125\delta}\mathscr{L}^3\ll y\mathscr{L}^{-A}+y\mathscr{L}^3\exp\left(-0.125\delta\mathscr{L}\right).
\end{align*}
Using the conditions $\lambda\le\left(2\pi k^2x^{k-2}y^2\right)^{-1}$ and $y=x^{\frac{\scriptstyle5}{\scriptstyle8}}\mathscr{L}^{1.5A+0.25b+18}$, we have
\begin{align*}
T_0=&\left(\frac xy+\lambda x^k\right)q^{\frac12}\mathscr{L}^{A+3}\le\left(\frac{x}y+\frac{x^2}{2\pi k^2y^2}\right)\mathscr{L}^{A+0.5b+3}\le \frac{x^2}{y^2}\mathscr{L}^{A+0.5b+3}<x,
\end{align*}
Using this inequality, we estimate from below the parameter $\delta=\delta(q,T_0)$:
\begin{align*}
\delta(q,T_0)=\frac{c_1}{\max\left(\ln q,\ \left(\ln(T_0+3)\ln\ln(T_0+3)\right)^\frac{3}{4}\right)}\ge \frac{c_1}{\max\left(b\ln\mathscr{L},\ \left(\mathscr{L}\ln\mathscr{L}\right)^\frac{3}{4}\right)}\ge c_1\mathscr{L}^{-0.76}.
\end{align*}
Therefore
\begin{align*}
\mathscr{A}_2&\ll y\mathscr{L}^{-A}+y\mathscr{L}^3\exp\left(-0.125c_1\mathscr{L}^{0.24}\right)\ll y\mathscr{L}^{-A}.
\end{align*}
Substituting this estimate and the estimate obtained for the sum $\mathscr{A}_1$ into (\ref{Formula W(alpha;x,y)<<A1+A2}), we obtain
\begin{align*}
|W(\alpha;x,y)|&\ll y\mathscr{L}^{-A}.
\end{align*}
From this estimate and (\ref{Formula otsenka summi W1}), in view of (\ref{formula Sk(alpha;x,y)=int+W1+W+R2}), we obtain the first assertion of the theorem.

\subsection{Proof of the theorem in the case $\lambda>\left(2\pi k^2x^{k-2}y^2\right)^{-1}$}
We introduce the parameter $H=2\pi k^2\lambda x^{k-1}y$, which for $\lambda>\left(2\pi k^2x^{k-2}y^2\right)^{-1}$ satisfies the condition
\begin{equation}\label{Formula H>x/y}
H\ge \frac{x}{y}.
\end{equation}
We divide all zeros $\rho=\beta+i\gamma$ of the function $L(s,\chi)$ with the condition $|\gamma|\le T_0$ into sets $D_1$, $D_2$ and $D_3$ as follows:
\begin{align*}
&D_1=\left\{\rho:\ -T_0\le\gamma<-2\pi k\lambda x^k-H\right\},  \\
&D_2=\left\{\rho:\ -2\pi k\lambda x^k-H\le\gamma\le-2\pi k\lambda(x-y)^k+H\right\}, \\
& D_3=\left\{\rho:\ -2\pi k\lambda(x-y)^k+H<\gamma\le T_0 \right\}.
\end{align*}
Denoting by $\EuScript{W}_j(\lambda,\chi)$, $j=1,\ 2,\ 3$, the sum of the moduli of the integral $I(\rho,\lambda)$ over the zeros $\rho$ belonging to the set $D_j$, we represent the sum $\EuScript{W}(\lambda,\chi)$ in the form
\begin{align}\label{Formula W(chi,lambda)=W1+W2+W3}
&\EuScript{W}(\lambda,\chi)=\EuScript{W}_1(\lambda,\chi)+\EuScript{W}_2(\lambda,\chi)+\EuScript{W}_3(\chi,\lambda).
\end{align}

{\bf Estimate of $\EuScript{W}_1(\lambda,\chi)$}. Adding the term $2\pi k\lambda u^k$, $x-y\leq u \le x$, to all three members of the inequality defining the set $D_1$, we obtain
\begin{align*}
&D_1=\left\{\rho:\ -T_0+2\pi k\lambda u^k\le\gamma+2\pi k\lambda u^k<-2\pi k\lambda x^k+2\pi k\lambda u^k-H\right\},
\end{align*}
On the interval $x-y\leq u\le x$ the function $2\pi k\lambda u^k$ is monotonically increasing, therefore for the right boundary of the set $D_1$ we have
$$
-2\pi k\lambda x^k+2\pi k\lambda u^k-H\le-H.
$$
Consequently, if $\rho$ belongs to the set $D_1$, then the inequality
$$
\gamma+2\pi k\lambda u^k<-H
$$
holds, therefore for the monotonically increasing function $\gamma+2\pi k\lambda u^k$ on the interval $x-y\le u\le x$ the relation
\begin{align*}
&\min|\gamma+2\pi k\lambda u^k|=-\max(\gamma+2\pi k\lambda u^k)=-\gamma-2\pi k\lambda x^k\ge H,
\quad \text{if}  \quad \rho\in D_1,
\end{align*}
holds. Hence, taking into account estimate (\ref{Formula otsenka I(rho,lambda)}), we find
$$
\EuScript{W}_1(\lambda,\chi)\ll\sum_{\rho \in D_1}\frac{x^\beta }{-\gamma-2\pi k\lambda x^k}.
$$
We divide all zeros in the set
$$
D_1=\left\{\rho:\ H\le-\gamma-2\pi k\lambda x^k<T_0-2\pi k\lambda x^k\right\},
$$
into classes $D_{11},\ \ldots,\ D_{1t}$, $t\ll T_0H^{-1}$, as follows: in the class $D_{1n}$ we place those zeros $\rho$ for which the conditions
\begin{align*}
nH<&-\gamma-2\pi k\lambda x^k\le(n+1)H,\hspace{26pt}\text{ if}\quad 1\le n\le t.
\end{align*}
hold. Therefore we have
\begin{align*}
\EuScript{W}_1(\lambda,\chi)\ll&\sum_{n=1}^t\sum_{\rho\in D_{1n}}\frac{x^\beta}{-\gamma-2\pi k\lambda x^k}\le\sum_{n=1}^t\sum_{\rho\in D_{1n}}\frac{x^\beta }{nH}\le\frac{\mathscr{L}}H\max_{1\le n\le t}\sum_{\rho\in D_{1n}}x^{\beta}\le\\
&\le\frac{\mathscr{L}}H\max_{|T|\le T_0}\sum_{T-H<\gamma\le T}x^\beta\ll\frac{\mathscr{L}}{\lambda x^{k-1}y}\max_{|T|\le T_0}\sum_{T-\lambda x^{k-1}y<\gamma\le T}x^\beta.
\end{align*}

{\bf Estimate of $\EuScript{W}_3(\lambda,\chi)$}. Adding the term $2\pi k\lambda u^k$, $x-y\leq u\le x$, to all three members of the inequality defining the set $D_3$, we obtain
\begin{align*}
& D_3=\left\{\rho:\ 2\pi k\lambda u^k-2\pi k\lambda(x-y)^k+H<\gamma+2\pi k\lambda u^k\le T_0+2\pi k\lambda u^k \right\}.
\end{align*}
On the interval $x-y\leq u\le x$ the function $2\pi k\lambda u^k$ is monotonically increasing, therefore for the left boundary of the set $D_3$ we have
$$
2\pi k\lambda u^k-2\pi k\lambda(x-y)^k+H\ge H.
$$
Consequently, if $\rho$ belongs to the set $D_3$, then the inequality
$$
\gamma+2\pi k\lambda u^k>H
$$
holds, therefore for the monotonically increasing function $\gamma+2\pi k\lambda u^k$ on the interval $x-y\le u\le x$ the relation
\begin{align*}
&\min|\gamma+2\pi k\lambda u^k|=\min(\gamma+2\pi k\lambda u^k)=\gamma+2\pi k\lambda(x-y)^k\ge H,
\quad \text{if}\quad\rho\in D_3.
\end{align*}
holds. We divide all zeros in the set
$$
 D_3=\left\{\rho:\ H<\gamma+2\pi k\lambda(x-y)^k\le T_0+2\pi k\lambda(x-y)^k\right\},
$$
into classes $D_{31},\ \ldots,\ D_{3t}$, $t\ll T_0H^{-1}$, as follows: in the class $D_{3n}$ we place those zeros $\rho$ for which the conditions
$$
nH<\gamma+2\pi k\lambda(x-y)^k\le(n+1)H,\hspace{26pt}\text{ if}\quad 1\le n\le t.
$$
hold. Therefore we have
\begin{align*}
\EuScript{W}_3(\lambda,\chi)\ll&\sum_{n=1}^t\sum_{\rho\in D_{3n}}\frac{x^\beta}{\gamma+2\pi k\lambda(x-y)^k}\le\sum_{n=1}^t\sum_{\rho\in D_{3n}}\frac{x^\beta }{nH}\le\frac{\mathscr{L}}H\max_{1\le n\le t}\sum_{\rho\in D_{3n}}x^{\beta}\le\\
&\le\frac{\mathscr{L}}H\max_{|T|\le T_0}\sum_{T-H<\gamma\le T}x^\beta\ll\frac{\mathscr{L}}{\lambda x^{k-1}y}\max_{|T|\le T_0}\sum_{T-\lambda x^{k-1}y<\gamma\le T}x^\beta.
\end{align*}

{\bf Estimate of $\EuScript{W}_2(\lambda,\chi)$}. Introducing the notation
$$
T_1=2\pi k\lambda x^k+H, \qquad H_1=2\pi k\lambda(x^k-(x-y)^k)+2H, \qquad T_2=T_1+2\pi k(k-1)\lambda u^k,
$$
where $x-y<u\le x$, we represent the set $D_2$ in the following forms:
\begin{align}\label{Formula mnozhestva D2}
D_2=\left\{\rho:\ T_1-H_1\le-\gamma\le T_1\right\}=\left\{\rho:\ T_2-H_1\le2\pi k(k-1)\lambda u^k-\gamma\le T_2\right\}.
\end{align}
Using the inequality
\begin{equation}\label{formula xk-(x-y)k<=kx(k-1)y}
(x-y)^k\ge x^k-kx^{k-1}y,
\end{equation}
which is a consequence of Lagrange's theorem on finite differences, and then the value of the parameter $H=2\pi k^2\lambda x^{k-1}y$, we have
\begin{align*}
T_2-H_1&=2\pi k(k-1)\lambda u^k+2\pi k\lambda(x-y)^k-H\ge2\pi k^2\lambda(x-y)^k-H\ge\\
&\ge2\pi k^2\lambda\left(x^k-kx^{k-1}y\right)-H=2\left(1-\frac{(k+1)y}{x}\right)\pi k^2\lambda x^k>\pi k^2\lambda x^k,
\end{align*}
Therefore, if $\rho\in D_2$, then using the third estimate (\ref{Formula otsenka I(rho,lambda)}), we find
\begin{align*}
|I(\rho,\lambda)|\ll\frac{x^\beta}{\sqrt{2\pi k(k-1)\lambda u^k-\gamma}}\le \frac{x^\beta}{\sqrt{T_2-H_1}}\le\frac{x^\beta}{\sqrt{\pi k^2\lambda x^k}}.
\end{align*}
Hence, and then from the representation of the set $D_2$ in the form (\ref{Formula mnozhestva D2}), we obtain
\begin{align}\label{Formula otsenka summi W2}
\EuScript{W}_2(\lambda,\chi)=&\sum_{\rho\in D_2}|I(\rho,\lambda)|\ll\sum_{T_2-H_1\le2\pi k(k-1)\lambda u^k-\gamma\le T_2}\frac{x^\beta}{\sqrt{\pi k^2\lambda x^k}}=
\frac1{\sqrt{\pi k^2\lambda x^k}}\sum_{T_1-H_1\le-\gamma\le T_1}x^\beta.
\end{align}
Using inequality (\ref{formula xk-(x-y)k<=kx(k-1)y}), we estimate from above $H_1$, the length of the set $D_2$, and from below its lower boundary $T_1-H_1$:
\begin{align*}
H_1&=2\pi k\lambda(x^k-(x-y)^k)+4\pi k^2\lambda x^{k-1}y\le6\pi k^2\lambda x^{k-1}y,\\
T_1-H_1&=2\pi k\lambda(x-y)^k-2\pi k^2\lambda x^{k-1}y\ge2\pi k\lambda x^k\left(1-\frac{2ky}x\right).
\end{align*}
Hence and from $T_1=2\pi k\lambda x^k\left(1+\frac{ky}x\right)$ it follows that in (\ref{Formula otsenka summi W2}) the boundaries of the sum over $\gamma$, i.e. $T_1-H_1$ and $T_1$, are quantities of order $\lambda x^k$, and its length $H_1$ is a quantity of order $\lambda x^{k-1}y$. We divide the summation interval $T_1-H_1\le\gamma\le T_1$ into at most $6\pi k^2$ intervals of the form
$$
\kappa\lambda x^k-\lambda x^{k-1}y<\gamma\le \kappa\lambda x^k,
$$
where the constant $\kappa$ takes a value from the interval $2\pi k\left(1-\frac{2ky}x\right)<\kappa\le2\pi k\left(1+\frac{ky}x\right)$, and represent the right-hand side of (\ref{Formula otsenka summi W2}) in the form
$$
\EuScript{W}_2(\lambda,\chi)\ll\frac{1}{\sqrt{\lambda x^k}}\sum_{\kappa\lambda x^k-\lambda x^{k-1}y<\gamma\le\kappa\lambda x^k}x^\beta.
$$
For convenience, without loss of generality, we shall assume that $\kappa=1$. Substituting this estimate, the estimates $\EuScript{W}_1(\lambda,\chi)$ and $\EuScript{W}_3(\lambda,\chi)$ into (\ref{Formula W(chi,lambda)=W1+W2+W3}) and bearing in mind that $\lambda x^k<T_0$, and then also $\lambda>\left(2\pi k^2x^{k-2}y^2\right)^{-1}$, we have
\begin{align*}
\EuScript{W}(\lambda,\chi)&\ll\frac{1}{\sqrt{\lambda qx^k}}\left(1+\frac{\mathscr{L}}{\sqrt{\lambda x^{k-2}}y}\right)\max_{|T|\le T_0}\sum_{T-\lambda x^{k-1}y<\gamma\le T}x^\beta\ll\frac{\mathscr{L}}{\sqrt{\lambda qx^k}}\max_{|T|\le T_0}\sum_{T-\lambda x^{k-1}y<\gamma\le T}x^\beta.
\end{align*}
Substituting this estimate into (\ref{Formula W(lambda,x,y)<<..sum W(lambda,chi)}), we obtain
$$
|W(\alpha;x,y)|\ll\frac{\mathscr{L}^2}{\sqrt{\lambda qx^k}}\max_{|T|\le T_0}\sum_{\chi\hspace{-7pt}\mod q}\sum_{T-\lambda x^{k-1}y<\gamma\le T}x^\beta=\frac{\mathscr{L}^2}{\sqrt{\lambda qx^k}}\max_{|T|\le T_0}\EuScript{V}_q(T,\lambda x^{k-1}y),
$$
where the sum $\EuScript{V}_q(T,U)$ is defined in (\ref{Formula W(alpha,x,y)<<frac()()...W(lambda)}), and in formula (\ref{Formula otsenka Vq(T,U)}) setting $U=\lambda x^{k-1}y$, we find
\begin{equation}\label{Formula otsenka summi W'}
|W(\alpha;x,y)|\ll\frac{\mathscr{L}^3}{\sqrt{\lambda qx^k}}\max_{|T|\le T_0}\max_{0.5\le u\le1-\delta}x^u\sum_{\chi\hspace{-7pt}\mod q}\left(N(u,T,\chi)-N(u,T-\lambda x^{k-1}y,\chi)\right),
\end{equation}
where $\delta=\delta(q,T_0)$. Further we shall assume that for the parameter $y$ the condition
\begin{equation}\label{Formula neravenstva y}
x^{1-\frac1{2k-1+\eta_k}}\mathscr{L}^{c_k}\le y\le x\mathscr{L}^{-A-0.5b-10}.
\end{equation}
holds. From the definition of the parameter $T_0$, the condition $\lambda>\left(2\pi k^2x^{k-2}y^2\right)^{-1}$ and (\ref{Formula neravenstva y}), we have
\begin{align*}
\frac{T}{(\lambda x^{k-1}y)^3}&\le\frac{T_0}{(\lambda x^{k-1}y)^3}
=\frac{\sqrt{q}\mathscr{L}^{A+3}}{\lambda^3x^{3k-4}y^4}+ \frac{\sqrt{q}\mathscr{L}^{A+3}}{\lambda^2x^{2k-3}y^3}\le\\
 &\le\frac{\left(2\pi k^2\right)^3y^2\mathscr{L}^{A+0.5b+3}}{x^2}+\frac{\left(2\pi k^2\right)^2y\mathscr{L}^{A+0.5b+3}}{x}\ll\mathscr{L}^{-7},
\end{align*}
it follows that in the sum over $\chi\hspace{-5pt}\mod q$ in (\ref{Formula otsenka summi W'}) the condition $\lambda x^{k-1}y\ge T^\frac{1}{3}$ holds, i.e. Lemma \ref{Lemma plotnostn-teor Zhan Tao} can be applied. Dividing the right-hand side of~\eqref{Formula otsenka summi W'} the interval $0.5\le u\le1-\delta$ into the intervals\linebreak $0.5\le u\le 0.75$ and $0.75\le u\le1-\delta$ we obtain two sums, which we shall denote by $\mathscr{B}_1$ and $\mathscr{B}_2$, respectively, and then setting in Lemma~\ref{Lemma plotnostn-teor Zhan Tao}
\begin{equation}\label{formula epsilon vib B2}
\varepsilon=\min\left(\frac{\eta_k}{8(k-1)},\ \frac{2k-5+\eta_k}2\delta\right),\qquad \delta=\delta(q,T_0),
\end{equation}
we have
\begin{align}
&|W(\alpha;x,y)|\ll\mathscr{B}_1+\mathscr{B}_2,\label{Formula otsenka summi B1}\\
&\mathscr{B}_1=\frac{\mathscr{L}^{12}}{\sqrt{\lambda qx^k}}\max_{0.5\le u\le0.75}x^u\left(q\lambda x^{k-1}y\right)^{\frac{4-4u}{3-2u}},\nonumber \\
&\mathscr{B}_2=\frac{\mathscr{L}^3}{\sqrt{q\lambda x^k}}\max_{0.75\le u\le1-\delta(q,T_0)}x^u\left(q\lambda x^{k-1}y\right)^{\frac{2}{u}(1-u)+\varepsilon}. \nonumber
\end{align}

\textbf{Estimate of $\mathscr{B}_1$}. We have
\begin{align*}
&\mathscr{B}_1= \frac{y^\frac12\mathscr{L}^{12}}{x^\frac12}\max_{0.5\le u\le0.75}\exp\left(f_1(u)\right),
&f_1(u)=u\ln x+\left(\frac{1}{u-1.5}+\frac32\right)\ln\left(q\lambda x^{k-1}y\right), \\
&f_1'(u)=\ln x-\frac{\ln\left(q\lambda x^{k-1}y\right)}{(u-1.5)^2},&f_1''(u)=\frac{2\ln\left(q\lambda x^{k-1}y\right)}{(u-1.5)^3}<0.\hspace{87pt}
\end{align*}
Let us show that $f_1'(0.5)>0$. Indeed, using the conditions $0\le\lambda\le\frac1{q\tau}$, $\tau=\frac{y^{2k-1}}{x^{k-1}\mathscr{L}^{b_1}}$ and $y\ge x^{1-\frac1{2k-1+\eta_k}}\mathscr{L}^{c_k}$, we have
\begin{align*}
-f_1'(0.5)&=\ln\left(q\lambda x^{k-2}y\right)\le\ln\left(\frac{x^{k-2}y}{\tau}\right)=\ln\left(\frac{x^{2k-3}\mathscr{L}^{b_1}}{y^{2k-2}}\right)\le  \ln\left(\frac{x^{2k-3}\mathscr{L}^{b_1}}{x^{2k-2-\frac{2k-2}{2k-1+\eta_k}}\mathscr{L}^{(2k-2)c_k}}\right)= \\
&=\ln\left(x^{-\frac{1+\eta_k}{2k-1+\eta_k}}\mathscr{L}^{b_1-(2k-2)c_k}\right)=-\frac{1+\eta_k}{2k-1+\eta_k}\mathscr{L}+(b_1-(2k-2)c_k)\ln\mathscr{L}<0.
\end{align*}
Further, to estimate $\mathscr{B}_1$ we consider two cases: $q\lambda x^{k-1}y<x^\frac9{16}$ and $q\lambda x^{k-1}y\ge x^\frac9{16}$.

\textbf{Case $q\lambda x^{k-1}y<x^\frac9{16}$}. We have
$$
f_1'(u)\ge f_1'(0.75)=\ln x-\frac{16}{9}\ln(q\lambda x^{k-1}y)>\ln x-\frac{16}{9}\ln x^\frac9{16}=0,
$$
i.e. $f_1'(u)$ is positive and $f_1(u)$ is an increasing function on the interval $0.5\le u\le0.75$. Using this property, and then relation (\ref{Formula neravenstva y}), we have
\begin{align*}
\mathscr{B}_1&=\frac{y^\frac12\mathscr{L}^{12}}{x^\frac12}\exp\left(f_1(0.75)\right)=\frac{y^\frac12\mathscr{L}^{12}}{x^\frac12}\exp\left(\frac{3}{4}\ln x+\frac16\ln(q\lambda x^{k-1}y)\right)=\left(q\lambda x^{k-1}y\right)^\frac16x^\frac14y^\frac12\mathscr{L}^{12}<\\
&<x^\frac{11}{32}y^\frac12\mathscr{L}^{12}= y\left(\frac{x^{1-\frac1{2k-1+\eta_k}}\mathscr{L}^{c_k}}{y}\right)^\frac12 x^{-\frac{10k+5\eta_k-21}{32(2k-1+\eta_k)}}\mathscr{L}^{12-0.5c_k}\ll y\mathscr{L}^{-A}.
\end{align*}

\textbf{Case $q\lambda x^{k-1}y\ge x^\frac9{16}$}. For convenience, introducing the notation
\begin{equation}\label{Formula znachenie y}
  y=x^{1-\mu^2}\mathscr{L}^c,
\end{equation}
we shall seek the largest value of the parameter $\mu$ and the smallest value of the parameter $c$ for which the estimate
\begin{equation}\label{Formula otsenka summi W'1-2sl}
\mathscr{B}_1\ll  y\mathscr{L}^{-A}.
\end{equation}
holds. Using the conditions $0\le\lambda\le\frac1{q\tau}$, $\tau=\frac{y^{2k-1}}{x^{k-1}\mathscr{L}^{b_1}}$ and $y\ge x^{1-\frac1{2k-1+\eta_k}}\mathscr{L}^{c_k}$, we have
$$
q\lambda x^{k-1}y\le\frac{x^{k-1}y}{\tau}=\left(\frac{x}{y}\right)^{2k-2}\mathscr{L}^{b_1}\le \left(x^{\frac1{2k-1+\eta_k}}\mathscr{L}^{-c_k}\right)^{2k-2}\mathscr{L}^{b_1}=x^{1-\frac{1+\eta_k}{2k-1+\eta_k}}\mathscr{L}^{b_1-(2k-2)c_k}.
$$
Consequently,
$$
f_1'(0.5)=\ln x-\ln(q\lambda x^{k-1}y)>0, \qquad  f_1'(0.75)=\ln x-\frac{16}{9}\ln(q\lambda x^{k-1}y)\le 0,
$$
and the point $u_0=\dfrac{3}{2}-\left(\dfrac{\ln(q\lambda x^{k-1}y)}{\ln x}\right)^\frac12$, where $f_1'(u_0)=0$, belongs to the interval of integration, and also on this interval the condition $ f''(u)<0$ holds. Consequently, on the interval $0.5\le u\le 0.75$ the graph of the function $f_1(u)$ is convex upward, therefore
\begin{align}
&\mathscr{B}_1=\frac{y^\frac12\mathscr{L}^{12}}{x^\frac12}\exp\left(f_1(u_0)\right)=
\frac{y^\frac12\mathscr{L}^{12}}{x^\frac12}\exp\left(u_0\ln x+\left(\frac{1}{u_0-1.5}+\frac32\right)\ln\left(q\lambda x^{k-1}y\right)\right)=\nonumber\\
&=\frac{y^\frac12\mathscr{L}^{12}}{x^\frac12}\exp\left(\left(\frac{3}{2}-\left(\frac{\ln(q\lambda x^{k-1}y)}{\ln x}\right)^\frac12\right)\ln x+\left(-\left(\frac{\ln x}{\ln(q\lambda x^{k-1}y)}\right)^\frac12+\frac32\right)\ln\left(q\lambda x^{k-1}y\right)\right)=\nonumber\\
&=xy^\frac12\mathscr{L}^{12}\exp\left(\frac32\ln\left(q\lambda x^{k-1}y\right)-2\left(\ln(q\lambda x^{k-1}y)\ln x\right)^\frac12\right)=xy^\frac12\mathscr{L}^{12}\exp\left(g\left(q\lambda x^{k-1}y\right)\right),\label{Formula otsenka summi B11}
\end{align}
where $g(t)=\frac32\ln t-2\left(\ln t\ln x\right)^\frac12$. From the conditions
$$
q\lambda x^{k-1}y\ge x^\frac9{16},\qquad \lambda\le\frac1{q\tau},\qquad \tau=\frac{y^{2k-1}}{x^{k-1}\mathscr{L}^{b_1}},
$$
and formula (\ref{Formula znachenie y}) it follows that
\begin{align}\label{Formula interval izm qlambdax..}
&x^\frac9{16}\le q\lambda x^{k-1}y\le\frac{x^{k-1}y}{\tau}=\frac{x^{2k-2}\mathscr{L}^{b_1}}{y^{2k-2}}=x^u\mathscr{L}^v, \qquad
u=(2k-2)\mu^2, \quad v=b_1-(2k-2)c.
\end{align}
Therefore, bearing in mind that $t=q\lambda x^{k-1}y\ge x^\frac9{16}$, we obtain
\begin{align*}
&g'(t)=\frac3{2t}-\frac{(\ln x)^\frac12}{t(\ln t)^\frac12}=\frac{9\ln t-4\ln x}{2t(\ln t)^\frac12(3(\ln t)^\frac12+2(\ln x)^\frac12)}>0,
\end{align*}
i.e. $g(t)$ is an increasing function on its interval of variation, therefore
$$
g\left(q\lambda x^{k-1}y\right)\le g\left(x^u\mathscr{L}^v\right).
$$
Hence, taking into account the relation
$$
y^{-\frac12}=x^{-\frac12+\frac{u}{4k-4}}\mathscr{L}^{-\frac{b_1-v}{4k-4}},
$$
we estimate the right-hand side of (\ref{Formula otsenka summi B11}) as follows:
\begin{align}
\mathscr{B}_1&\le xy^\frac12\mathscr{L}^{12}\exp\left(g\left(x^u\mathscr{L}^v\right)\right)=
yx^{\frac12+\frac{u}{4k-4}}\mathscr{L}^{12-\frac{b_1-v}{4k-4}}\exp\left(g\left(x^u\mathscr{L}^v\right)\right)=\nonumber\\
&=yx^{\frac12+\frac{u}{4k-4}}\mathscr{L}^{12-\frac{b_1-v}{4k-4}}\exp\left(\frac32\ln\left(x^u\mathscr{L}^v\right)- 2\left(\ln\left(x^u\mathscr{L}^v\right)\mathscr{L}\right)^\frac12\right)=\nonumber\\
&=y\mathscr{L}^{12-\frac{b_1-v}{4k-4}+\frac{3v}{2}}\cdot x^{\frac12+\frac{u}{4k-4}+\frac{3u}{2}}\exp\left(-2\left(u\mathscr{L}^2+v\mathscr{L}\ln\mathscr{L}\right)^\frac12\right)=\nonumber\\
&=y\mathscr{L}^{12-\frac{b_1}{4k-4}+\frac{v}{4k-4}+\frac{3v}{2}}x^{\frac{3u}2+\frac{u}{4k-4}+\frac12} \exp\left(-2\sqrt{u}\mathscr{L}\left(1+\frac{v\ln\mathscr{L}}{u\mathscr{L}}\right)^\frac12\right).\label{Formula otsenka summi W'1}
\end{align}
We use, for $|t|\le 0.1$, the formula
$$
(1+t)^\frac12=1+\frac{t}{2}+R(t),\qquad R(t,\theta)=-\frac{(1-\theta)t^2}{8(1+\theta t)^{\frac32}},
$$
which is obtained by expanding the function $f(t)=(1+t)^\frac12$ in a Maclaurin series with remainder term in Cauchy form
\begin{align*}
f(t)=1+&\sum_{n=1}^m\frac{\frac12\left(\frac12-1\right)\ldots\left(\frac12-(n-1)\right)}{n!}t^n+R_m(t,\theta),\qquad R_m(t,\theta)=\frac{f^{(m+1)}(\theta t)}{(m+1)!}(1-\theta)^mt^{m+1},\\
&0<\theta<1,\qquad f'(t)=\frac{1}{2}(1+t)^{-\frac12},\qquad f''(t)=-\frac{1}{4}(1+t)^{-\frac32},
\end{align*}
for $m=1$. We shall consider the remainder term $R_1(t,\theta)=-\dfrac{t^2(1-\theta)}{8(1+\theta tT)^\frac32}$ as a function of $\theta$. We have
\begin{align*}
&\frac{\partial R_1(t,\theta)}{\partial\theta}=\frac{((3-\theta)t+2)t^2}{16(1+\theta t)^\frac52}>0.
\end{align*}
Consequently,
$$
\min R(t,\theta)=R(t,0)=-\frac{t^2}{8}, \qquad (1+t)^\frac12=1+\frac{t}{2}+R(t)\ge 1+\frac{t}{2}-\frac{t^2}{8}.
$$
Therefore
\begin{align*}
&-2\sqrt{u}\mathscr{L}\left(1+\frac{v\ln\mathscr{L}}{u\mathscr{L}}\right)^\frac12\le -2\sqrt{u}\mathscr{L}\left(1+\frac{v\ln\mathscr{L}}{2u\mathscr{L}}-\frac{v^2\ln^2\mathscr{L}}{8u^2\mathscr{L}^2}\right)=\\
&=\ln x^{-2\sqrt{u}}+\ln\mathscr{L}^{-\frac{v}{\sqrt{u}}}+\frac{v^2\ln^2\mathscr{L}}{4u^\frac32\mathscr{L}}\le \ln x^{-2\sqrt{u}}+\ln\mathscr{L}^{-\frac{v}{\sqrt{u}}}+1.
\end{align*}
Substituting this estimate into the right-hand side of (\ref{Formula otsenka summi W'1}), we obtain
\begin{align}
\mathscr{B}_1\ll&y\mathscr{L}^{12-\frac{b_1}{4k-4}+\frac{v}{4k-4}+\frac{3v}{2}}x^{\frac{3u}2+\frac{u}{4k-4}+\frac12}\cdot x^{-2\sqrt{u}} \mathscr{L}^{-\frac{v}{\sqrt{u}}} =\nonumber\\
&=yx^{\left(\frac{3}{2}+\frac1{4k-4}\right)u-2\sqrt{u}+\frac12} \mathscr{L}^{12-\frac{b_1}{4k-4}+\left(\frac1{4k-4}+\frac32-\frac{1}{\sqrt{u}}\right)v}.\label{Formula otsenka summi W'1-2}
\end{align}
Using relation (\ref{Formula interval izm qlambdax..}), we express the exponents of $x$ and $\mathscr{L}$ on the right-hand side of the last inequality in terms of $\mu$ and $c$. Denoting these exponents by $\kappa(\mu)$ and $\omega(\mu,c)$, we have
\begin{align*}
\kappa(\mu)=&\left(\frac{3}{2}+\frac1{4k-4}\right)(2k-2)\mu^2-2\sqrt{(2k-2)}\ \mu+\frac12=\frac12\left((6k-5)\mu^2-4\sqrt{2k-2}\ \mu+1\right),\\
\omega(\mu,c)=&12-\frac{b_1}{4k-4}+\left(\frac1{4k-4}+\frac32-\frac{1}{\mu\sqrt{2k-2}}\right)(b_1-(2k-2)c)=\\
&=12+\left(\frac32-\frac{1}{\mu\sqrt{2k-2}}\right)b_1-\left(\frac{6k-5}{2}-\frac{\sqrt{2k-2}}{\mu}\right)c.
\end{align*}
The quantity $\kappa=\kappa(\mu)$ is a quadratic polynomial in $\mu$, with two positive roots, the larger of which is the number
\begin{align*}
\mu_k=\frac{2\sqrt{2k-2}+\sqrt{2k-3}}{6k-5}=\frac1{2\sqrt{2k-2}-\sqrt{2k-3}},
\end{align*}
whose square has the form
\begin{align}
\mu_k^2=&\frac{1}{2k-1+\eta_k}, &\eta_k=\frac{2}{4k-5+2\sqrt{(2k-2)(2k-3)}}. \label{Formula znachenie mu_k^2}
\end{align}
Consequently, the largest value of the parameter $\mu$ for which $\kappa=\kappa(\mu)$, the exponent of $x$ in estimate (\ref{Formula otsenka summi W'1-2}), equals zero is the number $\mu_k$, and $\omega(\mu,c)$, the exponent of $\mathscr{L}$ at $\mu=\mu_k$, takes the form
\begin{align*}
\omega(\mu_k,c)&=12+\left(\frac{\sqrt{2k-3}}{\sqrt{2k-2}}-\frac12\right)b_1-\left(\sqrt{(2k-2)(2k-3)}-\frac{2k-3}{2}\right)c,
\end{align*}
and for $c=c_k$, where
\begin{align*}
&c_k=\frac{2A+24+\left(\frac{2\sqrt{2k-3}}{\sqrt{2k-2}}-1\right)b_1}{2\sqrt{(2k-2)(2k-3)}-(2k-3)},
\end{align*}
the equality $\omega(\mu_k,c_k)=-A$ holds, i.e. estimate (\ref{Formula otsenka summi W'1-2}) turns into estimate (\ref{Formula otsenka summi W'1-2sl}). Hence, from the representations of the parameter $y$ in the form (\ref{Formula znachenie y}) and of the parameter $\mu_k^2$ in the form (\ref{Formula znachenie mu_k^2}) it follows that estimate (\ref{Formula otsenka summi W'1-2sl}) holds for
$$
y\ge x^{1-\frac{1}{2k-1+\eta_k}}\mathscr{L}^{c_k},  \qquad \eta_k=\frac{2}{4k-5+2\sqrt{(2k-2)(2k-3)}}.
$$

\textbf{Estimate of $\mathscr{B}_2$}. We have
\begin{align*}
&\mathscr{B}_2=\frac{y^\frac12\mathscr{L}^3}{x^\frac12}\left(q\lambda x^{k-1}y\right)^{-\frac12}\max_{0.75\le u\le1-\delta}f_2(u), \qquad f_2(u)=x^u\left(q\lambda x^{k-1}y\right)^{\frac{2}{u}(1-u)+\varepsilon},\qquad \delta=\delta(q,T_0).
\end{align*}
In $\mathscr{B}_2$ the function $f_2(u)$ and its second derivative are positive:
\begin{align*}
f_2'(u)=&f_2(u)\left(\ln x-\frac2{u^2}\ln\left(q\lambda x^{k-1}y\right)\right),\quad
f_2''(u)=f_2(u)\left(\left(\ln x-\frac2{u^2}\ln\left(q\lambda x^{k-1}y\right)\right)^2+\right.\\
&\left.+\frac4{u^3}\ln\left(q\lambda x^{k-1}y\right)\right)\ge\frac{4 f_2(u)}{u^3}\ln\left(q\lambda x^{k-1}y\right)\ge\frac{4 f_2(u)}{u^3}\ln\left(\frac{qx}{2\pi k^2y}\right)>0,
\end{align*}
i.e. the graph of the function $f_2(u)$ is convex downward, therefore
\begin{align*}
\mathscr{B}_2\le&\frac{y^\frac12\mathscr{L}^2}{x^\frac12}\left(q\lambda x^{k-1}y\right)^{-\frac12}\left(f_2(0.75)+f_2(1-\delta)\right)=\\
&=\frac{y^\frac12\mathscr{L}^2.}{x^\frac12}\left(x^\frac34\left(q\lambda x^{k-1}y\right)^{\frac16+\varepsilon}+ x^{1-\delta }\left(q\lambda x^{k-1}y\right)^{-\frac12+\frac{2\delta}{1-\delta}+\varepsilon}\right).
\end{align*}
Using the condition of the case under consideration
$$
\frac{1}{2\pi k^2x^{k-2}y^2}<\lambda\le\frac{1}{q\tau}=\frac{x^{k-1}\mathscr{L}^{b_1}}{qy^{2k-1}},
$$
we have
\begin{align*}
&\mathscr{B}_2\ll\frac{y^\frac12\mathscr{L}^2}{x^\frac12}\left(x^\frac34\left(\frac{x^{2k-2}\mathscr{L}^{b_1}}{y^{2k-2}}\right)^{\frac{\scriptstyle1}{\scriptstyle6}+\varepsilon}+ x^{1-\delta}\left(\frac{x}y\right)^{-\frac{\scriptstyle1}{\scriptstyle2}+\frac{\scriptstyle{2\delta}}{\scriptstyle{1-\delta}}+\varepsilon}\right) =y\mathscr{L}^2\cdot\Delta(k,\varepsilon)\\
&\Delta(k,\varepsilon)=\left(\frac{x^{2k-2+\frac{\scriptstyle3}{\scriptstyle{2(1+6\varepsilon)}}}\mathscr{L}^{b_1}}{y^{2k-2+ \frac{\scriptstyle3}{\scriptstyle{1+6\varepsilon}}}}\right)^{\frac{\scriptstyle1}{\scriptstyle6}+\varepsilon}+ \left(\frac{x^\frac{\scriptstyle{\delta+\delta^2+(1-\delta)\varepsilon}}{\scriptstyle{2\delta+(1-\delta)\varepsilon}}}{y} \right)^{\frac{\scriptstyle{2\delta}}{\scriptstyle{1-\delta}}+\varepsilon}=\\
&=\left(\frac{x^{1-\frac{\scriptstyle3}{\scriptstyle{4k+2+24(k-1)\varepsilon}}} \mathscr{L}^{\frac{\scriptstyle{b_1+6\varepsilon}}{\scriptstyle{2k+1+12(k-1)\varepsilon}}}}{y}\right)^{\frac{\scriptstyle{2k+1}}{\scriptstyle6} +2(k-1)\varepsilon}+\left(\frac{x^\frac{\scriptstyle{\delta+\delta^2+(1-\delta)\varepsilon}}{\scriptstyle{2\delta+(1-\delta)\varepsilon}}}{y} \right)^{\frac{\scriptstyle{2\delta}}{\scriptstyle{1-\delta}}+\varepsilon}=\\
&=\left(\frac{x^{1-\frac{\scriptstyle1}{\scriptstyle{2k-1+\eta_k}}}\mathscr{L}^{c_k}}yx^{g(k,\varepsilon)}\mathscr{L}^{h(k, \varepsilon)}\right)^{\frac{\scriptstyle{2k+1}}{\scriptstyle6}+2(k-1)\varepsilon} +\left(\frac{x^{1-\frac{\scriptstyle1}{\scriptstyle{2k-1+\eta_k}}}\mathscr{L}^{c_k}}y x^{f(\delta,\varepsilon)}\mathscr{L}^{-c_k}\right)^{\frac{\scriptstyle{2\delta}}{\scriptstyle{1-\delta}}+\varepsilon},
\end{align*}
where
\begin{align*}
&g(k,\varepsilon)=-\frac{2k-5+3\eta_k-24(k-1)\varepsilon}{2(2k-1+\eta_k)(2k+1+12(k-1)\varepsilon)},\\
&h(k,\varepsilon)=\frac{b_1+6\varepsilon}{2k+1+12(k-1)\varepsilon}-c_k,\\
&f(\delta,\varepsilon)=-\frac{\delta(1-\delta)}{2\delta+(1-\delta)\varepsilon}+\frac{1}{2k-1+\eta_3}.
\end{align*}
Hence, bearing in mind that $y\ge x^{1-\frac1{2k-1+\eta_k}}\mathscr{L}^{c_k}$ and $c_k>0$, we obtain
\begin{align}\label{Formula otsenka summi B2}
\Delta(k,\varepsilon)&\ll \left(x^{g(k,\varepsilon)}\mathscr{L}^{h(k, \varepsilon)}\right)^\frac{\scriptstyle{2k+1+12(k-1)\varepsilon}}{\scriptstyle6}+ x^{f(\delta,\varepsilon)\frac{\scriptstyle{2\delta+(1-\delta)\varepsilon}}{\scriptstyle{1-\delta}}}.
\end{align}
Using the explicit value of the parameter $\varepsilon$, i.e. formula (\ref{formula epsilon vib B2}), the conditions $\eta_k>0$, and $k\ge 3$, we find
\begin{align*}
g(k,\varepsilon)\frac{2k+1+12(k-1)\varepsilon}6&=-\frac{2k-5}{12(2k-1+\eta_k)}+\left(\varepsilon-\frac{\eta_k}{8(k-1)}\right) \frac{2k-2}{2k-1+\eta_k}\le\\
&\le-\frac{2k-5}{12(2k-1+\eta_k)}< -\frac{2k-5}{12(2k-1)}=-\frac1{12}+\frac1{6k-3}\le-\frac1{60},\\
f(\delta,\varepsilon)\frac{2\delta+(1-\delta)\varepsilon}{1-\delta}&=
-\delta+\frac{2\delta+(1-\delta)\varepsilon}{(1-\delta)(2k-1+\eta_k)}\le-\delta+\frac{2\delta+\varepsilon}{\frac56(2k-1+\eta_k)}=\\
&=-\frac25\delta+\frac6{5(2k-1+\eta_k)}\left(\varepsilon-\frac{2k-5+\eta_k}2\delta\right)\le-\frac25\delta.
\end{align*}
With the help of these two inequalities and the relation
\begin{align*}
h(k,\varepsilon)\left(\frac{2k+1}6+2(k-1)\varepsilon\right)=\frac{b_1-(2k+1)c_k+(6-12(k-1)c_k)\varepsilon}6
\end{align*}
we represent estimate (\ref{Formula otsenka summi B2}) in the form
\begin{align*}
\Delta(k,\varepsilon)&\ll x^{-\frac1{60}}\mathscr{L}^\frac{b_1-(2k+1)c_k+(6-2(k-1)c_k)\varepsilon}6+ x^{-0.4\delta}\ll \mathscr{L}^{-A-2}+\exp\left(-0.4\delta\mathscr{L}\right).
\end{align*}
Using the conditions
$$
\lambda\le\frac1{q\tau},\qquad \tau=\frac{y^{2k-1}}{x^{k-1}\mathscr{L}^{b_1}},\qquad y\ge x^{1-\frac1{2k-1+\eta_k}}\mathscr{L}^{c_k},
$$
we have
\begin{align*}
T_0\le\left(\frac{x\sqrt{q}}y+\frac{x^{2k-1}\mathscr{L}^{b_1}}{y^{2k-1}}\right)\mathscr{L}^{A+3}\ll\frac{x^{2k-1}}{y^{2k-1}}\mathscr{L}^{A+b_1+3}\le x^{\frac{2k-1}{2k-1+\eta_k}}\mathscr{L}^{A+b_1+3-(2k-1)c_k}<0.1x,
\end{align*}
Using this inequality, we estimate from below the parameter $\delta=\delta(q,T_0)$:
\begin{align*}
\delta(q,T_0)=\frac{c_1}{\max\left(\ln q,\ \left(\ln(T_0+3)\ln\ln(T_0+3)\right)^\frac{3}{4}\right)}\ge \frac{c_1}{\max\left(b\ln\mathscr{L},\ \left(\mathscr{L}\ln\mathscr{L}\right)^\frac{3}{4}\right)}\ge c_1\mathscr{L}^{-0.76}.
\end{align*}
Therefore
\begin{align*}
\mathscr{B}_2\ll y\mathscr{L}^2\cdot\Delta(k,\varepsilon)\ll y\mathscr{L}^2\left(\mathscr{L}^{-A-2}+\exp\left(-0.4c_1\mathscr{L}^{0.24}\right)\right)\ll y\mathscr{L}^{-A}.
\end{align*}
Substituting this estimate and the estimate for $\mathscr{B}_1$ into (\ref{Formula otsenka summi B1}), we obtain
\begin{align}
|W(\alpha;x,y)|&\ll y\mathscr{L}^{-A}.\label{Formula otsenka W(lambda,x,y)>}
\end{align}
Using estimate (\ref{Formula otsenka tau(chi,a,k)}), Lemma \ref{Lemma otsenka trig integ po pervoy proizv} on the estimate of the integral by the magnitude of the first derivative, the conditions $\lambda>\left(2\pi k^2x^{k-2}y^2\right)^{-1}$ and (\ref{Formula neravenstva y}), we have
$$
\left|\frac{\tau(\chi_0,a,k)}{\varphi(q)}\int_{x-y}^xe(\lambda u^k)du\right|\ll\frac{\mathscr{L}}{\sqrt{q}}\cdot\frac{1}{\lambda x^{k-1}}\ll
\frac{y^2\mathscr{L}}{\sqrt{q}x}\ll y\mathscr{L}^{-A}.
$$
Substituting this estimate, as well as estimates (\ref{Formula otsenka W(lambda,x,y)>}) and (\ref{Formula otsenka summi W1}) into (\ref{formula Sk(alpha;x,y)=int+W1+W+R2}), we obtain the second assertion of the theorem.

\end{document}